\documentclass{amsart}

\usepackage{fancyhdr}
\usepackage{lastpage}
\usepackage{stmaryrd,yhmath}

\newcommand{\bdis}{\begin{displaymath}}
\newcommand{\edis}{\end{displaymath}}
\newcommand{\be}{\begin{equation}}
\newcommand{\ee}{\end{equation}}
\newcommand{\mbb}{\mathbb}
\newcommand{\mcal}{\mathcal}

\newcommand{\vp}{\varphi}

\newcommand{\vth}{\vartheta}

\newcommand{\zf}{\zeta\left(\frac{1}{2}+it\right)}

\newtheorem{lemma}[]{Lemma}

\theoremstyle{definition}

\newtheorem{cor}[]{Corollary}

\theoremstyle{remark}
\newtheorem{remark}[]{Remark}

\newtheorem*{mydef1}{{\bf Theorem}}

\newtheorem*{mydef3}{{\bf Problem}}

\numberwithin{equation}{section}

\begin{document}

\title{Jacob's ladders and the $\tilde{Z}^2$-transformation of a polynomials in $\ln \vp_1(t)$}

\author{Jan Moser}

\address{Department of Mathematical Analysis and Numerical Mathematics, Comenius University, Mlynska Dolina M105, 842 48 Bratislava, SLOVAKIA}

\email{jan.mozer@fmph.uniba.sk}

\keywords{Riemann zeta-function}

\begin{abstract}
It is proved in this paper that there is a nonlocal asymptotic splitting (in the integral sense) of the function $Z^4(t)$ into two factors. The
corresponding formula cannot be obtained in the known theories of Balasubramanian, Heath-Brown and Ivic.
\end{abstract}

\maketitle

\section{Resluts: nonlocal splitting of the function $Z^4(t)$}

\subsection{}

Let us remind that Hardy and Littlewood started in 1926 to study the following integral
\be \label{1.1}
\int_1^T\left| \zf\right|^4{\rm d}t=\int_1^T Z^4(t){\rm d}t,
\ee
where
\bdis
Z(t)=e^{i\vth(t)}\zf,\ \vth(t)=-\frac{t}{2}\ln\pi+\text{Im}\ln\Gamma\left(\frac{1}{4}+i\frac{t}{2}\right),
\edis
and they derived the following estimate (see \cite{1}, pp. 41, 59, \cite{16}, p. 124)
\be \label{1.2}
\int_1^T \left|\zf\right|^4{\rm d}t=\mcal{O}(T\ln^4 T) .
\ee
In 1926 Ingham derived the asymptotic formula
\be \label{1.3}
\int_1^T \left|\zf\right|^4{\rm d}t=\frac{1}{2\pi^2}T\ln^4 T+\mcal{O}(T\ln^3 T)
\ee
(see \cite{2}, p. 277, \cite{16}, p. 125). Let us remind, finally, the Ingham - Heath-Brown formula (see \cite{3}, p.129)
\be \label{1.4}
\int_0^T Z^4(t){\rm d}t=T\sum_{k=0}^4 C_k\ln^{4-k} T+\mcal{O}(T^{7/8+\epsilon}),\ C_0=\frac{1}{2\pi^2}
\ee
which improved the Ingham formula (\ref{1.3}).

\subsection{}

From (\ref{1.4}) the formula
\be \label{1.5}
\int_T^{T+U_0} Z^4(t){\rm d}t=\left. t\sum_{k=0}^4 C_k\ln^{4-k} t\right|_T^{T+U_0}+\mcal{O}(T^{7/8+\epsilon})
\ee
follows, where $U_0=T^{7/8+2\epsilon}$.

\begin{mydef3}
Which is the behaviour of the integral in (\ref{1.5}) under the translation (in the asymptotic sense)
\be \label{1.6}
\mcal{T}^{(-)}:\ [T,T+U_0]\mapsto [\vp_1(T),\vp_1(T+U_0)]
\ee
where $\vp_1(T),\ T\geq T_0[\vp_1]$ is the Jacob's ladder.
\end{mydef3}

The following theorem gives us the answer.

\begin{mydef1}
Let
\bdis
\tilde{Z}^2(t)=\frac{Z^2(t)}{2\Phi^\prime_\vp[\vp[t]]}
\edis
(see \cite{14}, (5.1), (5.2)). Then
\begin{eqnarray} \label{1.7}
& &
\int_{\vp_1(T)}^{\vp_1(T+U_0)}Z^4(t){\rm d}t= \\
& &
=\int_T^{T+U_0}\left\{\sum_{k=0}^4 D_k\ln^{4-k}\vp_1(t)\right\}\tilde{Z}^2(t){\rm d}t+\mcal{O}(T^{7/8+\epsilon}), \nonumber
\end{eqnarray}
where
\be \label{1.8}
(D_0,\dots ,D_4)=(C_0,C_1+4C_0,C_2+3C_1,C_3+2C_2,C_4+C_3) ,
\ee
and
\be \label{A}
\vp_1(T+U_0)-\vp_1(T)\sim U_0,\ \vp_1(T+U_0)<T; \ (^-) \tag{A}
\ee
\be \label{B}
\rho\{ [T,T+U_0];[\vp_1(T),\vp_1(T+U_0)]\}\sim (1-c)\pi(T)\to\infty \tag{B}
\ee
if $T\to\infty$ and $\rho$ denotes the distance of the corresponding segments, $c$ is the Euler constant and $\pi(T)$ is the prime-counting
function.
\end{mydef1}

Since
\bdis
\int_{\vp_1(T)}^{\vp_1(T+U_0)}Z^4(t){\rm d}t>AU_0\ln^4(T)=AT^{7/8+2\epsilon}\ln^4 T,
\edis
we obtain

\begin{cor}
\be \label{1.9}
\int_{\vp_1(T)}^{\vp_1(T+U_0)}Z^4(t){\rm d}t\sim \int_T^{T+U_0}\left\{\sum_{k=0}^4 D_k\ln^{4-k}\vp_1(t)\right\}\tilde{Z}^2(t){\rm d}t .
\ee
\end{cor}
Next we have

\begin{cor}
Let $[T,T+U_0]=\vp_1\{[\mathring{T},\widering{T+U_0}]\}$. Then by the translation
\bdis
\mcal{T}^{(+)}:\ [T,T+U_0]\mapsto [\mathring{T},\widering{T+U_0}]
\edis
the following splitting
\be \label{1.10}
\int_T^{T+U_0}Z^4(t){\rm d}t\sim \int_{\mathring{T}}^{\widering{T+U_0}}\left\{\sum_{k=0}^4 D_k\ln^{4-k}\vp_1(t)\right\}\tilde{Z}^2(t){\rm d}t
\ee
where
\be \label{A1}
T+U_0<\mathring{T}; \ (^+) \tag{A1}
\ee
\be\label{B1}
\rho\{[T,T+U_0];[\mathring{T},\widering{T+U_0}]\}\sim (1-c)\pi(T) \tag{B1}
\ee
is generated.
\end{cor}

\begin{remark}
Due to the formulae (1.7), (1.9), (1.10) the following effect - the nonlocal asymptotic splitting (in the integral sense) of the function $Z^4(t)$ into two
factors - is expressed.
\end{remark}

\begin{remark}
It is obvious that the formulae (1.7), (1.9) and (1.10) remain valid for all
\bdis
U_0=T^{\omega+2\epsilon},\ \omega<\frac{7}{8} ,
\edis
where $\omega$ is an arbitrary improvement of the exponent $7/8$ which will be shown (since the integral on the right-hand side of (1.9) holds true
for all $U_0\in (0,T/\ln T]$).
\end{remark}

\begin{remark}
It is clear that the formulae (1.7), (1.9) and (1.10) cannot be obtained in known theories of Balasubramanian, Heath-Brown and Ivic (comp. \cite{3}).
\end{remark}

This paper is a continuation of the paper series \cite{5}-\cite{15}.

\section{Universal character of $|\zf|^2$}

Let us remind that
\bdis
\tilde{Z}^2(t)=\frac{{\rm d}\vp_1(t)}{{\rm d}t},\ \vp_1(t)=\frac{1}{2}\vp(t) ,
\edis
where
\be \label{2.1}
\tilde{Z}^2(t)=\frac{Z^2(t)}{2\Phi^\prime_\vp[\vp(t)]}=\frac{Z^2(t)}{\left\{1+\mcal{O}\left(\frac{\ln\ln t}{\ln t}\right)\right\}\ln t},
\ee
(see \cite{15}, (5.1)-(5.3)). The following lemma holds true.

\begin{lemma}
For every integrable function (in the Lebesgue sense) $f(x),\ x\in[\vp_1(T),\vp_1(T+U)]$ the following is true
\be \label{2.2}
\int_T^{T+U} f[\vp_1(t)]\tilde{Z}^2(t){\rm d}t=\int_{\vp_1(T)}^{\vp_1(T+U)}f(x){\rm d}x,\ U\in (0,T/\ln T] ,
\ee
where $t-\vp_1(t)\sim(1-c)\pi(t)$.
\end{lemma}

\begin{remark}
The formula (2.2) is true also in the case when the integral on the right-hand side of (2.2) is convergent but non-absolutely (in the
Riemann sense).
\end{remark}

If $\vp_1\{ [\mathring{T},\widering{T+U}]\}=[T,T+U]$ then we have the following formula (see (2.2)).

\begin{lemma}
\be \label{2.3}
\int_{\mathring{T}}^{\widering{T+U}} f[\vp_1(t)]\tilde{Z}^2(t){\rm d}t=\int_T^{T+U} f(x){\rm d}x,\ U\in (0,T/\ln T] .
\ee
\end{lemma}

\begin{remark}
By the formula (2.2) the function
\bdis
\left|\zf\right|^2=Z^2(t)
\edis
is connected with all \emph{normal} functions of the Analysis (on the other cases of the universality of $\zeta(s)$ see \cite{4}, pp. 130-135).
\end{remark}

\section{On $\tilde{Z}^2$-transformation of some polynomials in $\ln\vp_1(t)$}

The following lemma is true.

\begin{lemma}
\begin{eqnarray} \label{3.1}
& &
\int_T^{T+U}\left\{\sum_{k=0}^4 A_k\ln^{4-k}\vp_1(t)\right\}\tilde{Z}^2(t){\rm d}t= \\
& &
= x\left. \sum_{k=0}^4 B_k\ln^{4-k}x\right|_{x=\vp_1(t)}^{x=\vp_1(T+U)}, \ U\in (0,T/\ln T] , \nonumber
\end{eqnarray}
where
\begin{eqnarray} \label{3.2}
B_0 & = & A_0 ,\nonumber  \\
B_1 & = & -4A_0+A_1, \nonumber  \\
B_2 & = & 12A_0-3A_1+A_2, \\
B_3 & = & -24A_0+6A_1-2A_2+A_3, \nonumber \\
B_4 & = & 24A_0-6A_1+2A_2-A_3+A_4 , \nonumber
\end{eqnarray}
i.e.
\be \label{3.3}
(A_0,\dots ,A_4)=(B_0,B_1+4B_0,B_2+3B_1,B_3+2B_2,B_4+B_3) .
\ee
\end{lemma}

\begin{proof}
The expressions (3.1), (3.2) follows from the formulae (see (2.2)):
\begin{eqnarray*}
& &
\int_T^{T+U}\tilde{Z}^2(t){\rm d}t=\left. x\right|_{\vp_1(T)}^{\vp_1(T+U)}, \\
& &
\int_T^{T+U}\tilde{Z}^2(t)\ln\vp_1(t){\rm d}t=\left. x (\ln x-1)\right|_{\vp_1(T)}^{\vp_1(T+U)} , \\
& &
\int_T^{T+U}\tilde{Z}^2(t)\ln^2\vp_1(t){\rm d}t=\left. x (\ln^2 x-2\ln x + 1)\right|_{\vp_1(T)}^{\vp_1(T+U)} , \\
& &
\int_T^{T+U}\tilde{Z}^2(t)\ln^3\vp_1(t){\rm d}t=\left. x (\ln^3 x-3\ln^2x+6\ln x-6)\right|_{\vp_1(T)}^{\vp_1(T+U)} , \\
& &
\int_T^{T+U}\tilde{Z}^2(t)\ln^4\vp_1(t){\rm d}t=\left. x (\ln^4 x-4\ln^3x+12\ln^2x-24\ln x+24)\right|_{\vp_1(T)}^{\vp_1(T+U)} ,
\end{eqnarray*}
\end{proof}

\section{Proof of the Theorem}

From (1.4) we obtain
\begin{eqnarray} \label{4.1}
& &
\int_{\vp_1(T)}^{\vp_1(T+U)} Z^4(t){\rm d}t= \\
& &
=t\left. \sum_{k=0}^4 C_k\ln^{4-k} t\right|_{\vp_1(T)}^{\vp_1(T+U)}+\mcal{O}(T^{7/8+\epsilon}) . \nonumber
\end{eqnarray}
Next, from (3.1)-(3.3), in the case
\bdis
(B_0,\dots ,B_4)=(C_0,\dots ,C_4)
\edis
we obtain (see (3.3), $B\to C,\ A\to D$)
\begin{eqnarray} \label{4.2}
& &
\int_T^{T+U_0}\left\{ \sum_{k=0}^4 D_k\ln^{4-k}\vp_1(t)\right\}\tilde{Z}^2(t){\rm d}t= \\
& &
=t\left. \sum_{k=0}^4 C_k\ln^{4-k} t\right|_{\vp_1(T)}^{\vp_1(T+U_0)} .
\end{eqnarray}
Then from (4.1), (4.2) we obtain (1.7). The expressions (\ref{A}), (\ref{B}) we obtain similarly to (C), (D), $k=1$ of Theorem in
\cite{11}.

\thanks{I would like to thank Michal Demetrian for helping me with the electronic version of this work.}

\end{document}